\documentclass[12pt]{article}
\usepackage{fullpage,graphicx,psfrag,amsmath,verbatim,xcolor}
\usepackage{tabu}
\usepackage[small,bf]{caption}
\usepackage[noend]{algpseudocode}
\usepackage{algorithm}
\usepackage[font=small,labelfont=bf,tableposition=top]{caption}
\usepackage{url}
\usepackage{wrapfig}
\usepackage{bbm}
\usepackage{multirow}
\DeclareCaptionLabelFormat{andtable}{#1~#2  \&  \tablename~\thetable}
\usepackage{tikz}
\usetikzlibrary{positioning,calc,shapes,arrows}
\tikzstyle{block} = [rectangle, draw,
text width=7em, text centered, minimum height=5em]
\tikzstyle{line} = [draw, -latex',line width=1mm]
\newcommand{\ones}{\mathbf 1}

\newcommand{\reals}{{\mbox{\bf R}}}

\newcommand{\BEAS}{\begin{eqnarray*}}
\newcommand{\EEAS}{\end{eqnarray*}}
\newcommand{\BEA}{\begin{eqnarray}}
\newcommand{\EEA}{\end{eqnarray}}
\newcommand{\BEQ}{\begin{equation}}
\newcommand{\EEQ}{\end{equation}}
\newcommand{\BIT}{\begin{itemize}}
\newcommand{\EIT}{\end{itemize}}

\newcommand{\Tr}{\mathop{\bf Tr}}

\newcommand{\prox}{\mathbf{prox}}

\newcommand{\argmin}{\mathop{\rm argmin}}

\newcommand{\eg}{{\it e.g.}}
\newcommand{\ie}{{\it i.e.}}

\newcounter{algorithmctr}[section]
\renewcommand{\thealgorithmctr}{\thesection.\arabic{algorithmctr}}
\newenvironment{algdesc}%
   {\refstepcounter{algorithmctr}
   \begin{list}{}{%
       \setlength{\rightmargin}{0\linewidth}%
       \setlength{\leftmargin}{.05\linewidth}}%
       \rmfamily\small
       \item[]{\setlength{\parskip}{0ex}\hrulefill\par%
        \nopagebreak{\bfseries\textsf{Algorithm \thealgorithmctr~}}}}%
   {{\setlength{\parskip}{-1ex}\nopagebreak\par\hrulefill} 
   \end{list}}

   {\refstepcounter{algorithmctr}
   \begin{list}{}{%
       \setlength{\rightmargin}{0\linewidth}%
       \setlength{\leftmargin}{.05\linewidth}}%
       \rmfamily\small
       \item[]{\setlength{\parskip}{0ex}\hrulefill\par%
        \nopagebreak{\bfseries\textsf{Algorithm}}}}%
   {{\setlength{\parskip}{-1ex}\nopagebreak\par\hrulefill} 
   \end{list}}

\title{Solving Large Multicommodity Network Flow Problems on GPUs}
\author{Fangzhao Zhang 
\and Stephen Boyd}
\begin{document}
\maketitle
\begin{abstract}
We consider the all-pairs multicommodity network flow problem on a 
network with capacitated edges.  The usual treatment keeps track of a 
separate flow for each source-destination pair on each edge; 
we rely on a more efficient formulation
in which flows with the same destination are aggregated, 
reducing the number of variables by a factor equal to the size of the 
network. Problems with hundreds of nodes, with
a total number of variables on the order of a million,
can be solved using standard generic interior-point methods on CPUs;
we focus on GPU-compatible algorithms that can solve such problems
much faster, and in addition scale to much larger problems, with 
up to a billion variables.
Our method relies on the primal-dual hybrid gradient
algorithm, and exploits several specific features of the problem for 
efficient GPU computation.
Numerical experiments show that our
primal-dual multicommodity network flow method accelerates
state of the art generic commercial solvers by 
$100\times$ to $1000\times$, and scales to problems that are much larger. 
We provide an open source implementation of our method.
\end{abstract}

\clearpage
\section{Multicommodity network flow optimization} \label{s-prob}
\subsection{Multicommodity network flow problem}\label{prob}
Our formulation of the multicommodity network flow (MCF) problem,
given below,
follows \cite{yin2019networkoptimizationunifiedpacket}.

\paragraph{Network.}
We consider a directed network with $n$ nodes and $m$ edges which is 
completely connected, \ie, there is a directed path between each 
pair of nodes. 
Let $A\in \reals^{n \times m}$ denote its incidence matrix, \ie,
\[
A_{i\ell} = \left\{ \begin{array}{rl} +1 &
\mbox{edge $\ell$ enters node $i$}\\
-1 &
\mbox{edge $\ell$ leaves node $i$}\\
0 & \mbox{otherwise}.
\end{array} \right.
\]
Edge $\ell$ has a positive capacity $c_\ell$.
The total flow on edge $\ell$ (to be defined below) cannot exceed $c_\ell$.

\paragraph{Traffic matrix.}
We consider the all-pairs multicommodity flow setting, \ie, 
there is traffic that originates at every node, destined for every other node.
We characterize the traffic between all source-destination pairs
via the traffic matrix $T\in \reals^{n\times n}$.
For any pair of distinct nodes $i,j$, $T_{ij}\geq 0$ 
is the traffic from (source) node $j$ to (destination) node $i$.  
There is no traffic from a node to itself; for mathematical convenience 
we define the diagonal traffic matrix entries as
$T_{ii}=-\sum_{j\neq i}T_{ij}$, the negative of the total
traffic with destination node $i$.
With this definition of the diagonal entries we have $T\ones = 0$,
where $\ones$ is the vector with all entries one.
 
\paragraph{Network utility.}
Let $u_{ij}$ denote the strictly concave increasing utility function 
for traffic from node $j$ to node $i$, for $j\neq i$.
We will assume utility functions are differentiable with 
domains $\reals_{++}$, the set of positive numbers.
(The methods we describe are readily extended to nondifferentiable 
utilities using subgradients instead of gradients.)
The total utility, which we wish to maximize, is 
$\sum_{i \neq j} u_{ij}(T_{ij})$.
For simplicity we take $u_{ii}=0$, so we can write the total
utility as 
\[
U(T) = \sum_{i,j} u_{ij}(T_{ij}).
\]
The domain of $U$ is 
$\mathcal T = \{T \mid T_{ij} >0~\mbox{for}~i\neq j\}$,
\ie, the traffic matrix must have positive off-diagonal entries.

Common examples of utility functions include the weighted log utility
$u(s)=w \log s $, and the weighted power utility
$u(s)=ws^\gamma$, with $\gamma \in (0,1)$,
where $w>0$ is the weight. 

\paragraph{Destination-based flow matrix.}
Following \cite{yin2019networkoptimizationunifiedpacket} we 
aggregate all flows with the same destination, considering it to be one 
commodity that is conserved at all nodes except the source and 
destination, but can be split and combined.
The commodity flows are given by the (destination-based) 
flow matrix $F\in\reals^{n\times m}$, where
$F_{i\ell}\geq 0$ denotes the flow on edge $\ell$ that is 
destined to node $i$. 
The edge capacity constraint can be expressed as 
$F^T\ones \leq c$, where the inequality is elementwise. 
A similar flow aggregation formulation, though source-based, was considered 
in \cite{bar2003origin}.

\paragraph{Flow conservation.}
The flow destined for node $i$ is conserved at all nodes $j \neq i$,
including the additional injection of traffic $T_{ij}$ that originates 
at node $j$ and is destined for node $i$.  This means that
\[
T_{ij}+\sum_\ell A_{j\ell}F_{i\ell}=0, \quad i,j=1,\ldots,n,\quad j\neq i.
\]
At the destination node, all traffic exits and we have
(using our definition of $T_{ii}$)
\[
T_{ii}+\sum_\ell A_{i\ell}F_{i\ell}=0, \quad i=1,\ldots,n.
\]
Combining these two, and using our specific definition of $T_{ii}$,
flow conservation can be compactly written in matrix notation as
\[
T+FA^T=0. 
\]

\paragraph{Multicommodity flow problem.}
In the MCF problem we seek a flow matrix 
that maximizes total network utility, subject to the 
edge capacity and flow conservation constraints.
This can be expressed as the problem
\BEQ\label{e-prob}
\begin{array}{ll}
\mbox{maximize} & U(T)\\
\mbox{subject to} & F\geq 0, \quad F^T\ones \leq c, \quad
T+FA^T=0,\\
\end{array}
\EEQ
with variables $F$ and $T$, and implicit constraint $T\in \mathcal T$.
The problem data are the network topology $A$, edge capacities $c$, 
and the traffic utility functions $u_{ij}$.

We can eliminate the traffic matrix $T$ using $T=-FA^T$ and state
the MCF problem in terms of the variable $F$ alone as
\BEQ\label{e-compact-prob}
\begin{array}{ll}
\mbox{maximize} & U(-FA^T) \\
\mbox{subject to} & F\geq 0, \quad F^T\ones \leq c,
\end{array}
\EEQ
with variable $F$, and implicit constraint $-FA^T \in \mathcal T$.
The number of scalar variables in this problem is $nm$.
For future use we define the feasible flow set as
\[
\mathcal F= \{F \mid F \geq 0,~F^T\ones \leq c \}.
\]

\paragraph{Existence and uniqueness of solution.}

First let us show the MCF problem \eqref{e-prob} is always feasible.
Consider a unit flow from each source to each destination,
over the shortest path, \ie, smallest number of edges, which exists since
the graph is completely connected. 
We denote this flow matrix as $F^\text{sp}$.
Now take $F=\alpha F^\text{sp}$, where 
$\alpha = 1/ \max_\ell(({F^\text{sp}}^T\ones)_\ell/c_\ell)>0$,
so we have $F^T\ones \leq c$. Evidently $F$ is feasible, 
and we have $T_{ij}= \alpha>0$ for $i \neq j$, so $T=-FA^T \in \mathcal T$.
This shows that the problem is always feasible.
Let $U^\text{sp}$ denote the corresponding objective function.

We can add the constraint $U(T)\geq U^\text{sp}$ to the problem,
without changing the solution set. 
With this addition, the feasible set is compact.
It follows that the MCF problem \eqref{e-prob} always has solution.  
The solution need not be unique.
The optimal $T$, however, is unique.
We also note that the argument above tells us that the implicit constraint
$T= -FA^T \in \mathcal T$ is redundant.

\paragraph{Solving MCF.}
The multicommodity flow problem \eqref{e-compact-prob} is convex 
\cite{boyd2004convex},
and so in principal can be efficiently solved.
In \cite{yin2019networkoptimizationunifiedpacket} the authors 
use standard generic interior-point solvers such as the commercial 
solver MOSEK \cite{mosek}, together with CVXPY \cite{diamond2016cvxpy},
to solve instances of the problem with
tens of nodes, and thousands of variables, in a few seconds on a CPU.
In this paper we introduce an algorithm for solving the MCF problem
that fully exploits GPUs.  For small and medium size problems
our method gives a substantial speedup over generic methods;
in addition it scales to much larger problems that cannot be solved 
by generic methods.

\subsection{Optimality condition and residual}
\paragraph{Optimality condition.} Let $\tilde {\mathcal F}$ denote 
the closure of the feasible set, including the implicit constraint 
$T=-FA^T \in \mathcal T$,
\[
\tilde {\mathcal F} = \mathcal F \cap 
\{F\mid -FA^T\in\textbf{cl}(\mathcal T)\},
\]
where $\textbf{cl}(\mathcal T)$ denotes the closure of $\mathcal T$.

Then $F$ is optimal for \eqref{e-compact-prob} if and only if 
$F\in \mathcal F$, $-FA^T \in \mathcal T$, and 
\[
\Tr (Z-F)^T G \geq 0
\]
holds for all $Z\in \tilde {\mathcal F}$,
where $G= \nabla_F (-U)(-FA^T)$
(see, \eg, \cite[\S 4.2.3]{boyd2004convex}).
We have $G=U'A$, where $U'_{ij}=u_{ij}'((-FA^T)_{ij})$. 

\paragraph{Optimality condition via projection onto $\mathcal F$.} For 
future use, we express the above optimality condition 
in terms of projection of a matrix $Q$ onto $\mathcal F$.
Let $\Pi$ denote Euclidean projection onto $\mathcal F$. 
Suppose $Q\in \reals^{n \times m}$, and set $F=\Pi (Q)$, so $F \in \mathcal F$.
Suppose in addition that $-FA^T \in \mathcal T$, so that 
$G=\nabla_F((-U)(-FA^T)$ exists.
Then $F$ is also Euclidean projection of $Q$ onto $\tilde {\mathcal F}$.
It follows that $\Tr(Z-F)^TG \geq 0$ for all $Z \in \tilde {\mathcal F}$,
so the optimality condition above holds, and $F$ is optimal.
Evidently it would hold if the weaker condition
\[
G=\gamma(F-Q)~\mbox{for some}~\gamma \geq 0
\]
holds.

Summarizing: $F$ is optimal if $F=\Pi(Q)$ for some $Q$, $-FA^T \in \mathcal T$,
and $G=\gamma(F-Q)$ for some $\gamma \geq 0$.
The converse is also true: If $F$ is optimal then 
$F=\Pi(Q)$ for some $Q$ with $-FA^T \in \mathcal T$
and $G=\gamma(F-Q)$ for some $\gamma \geq 0$.
(Indeed, this holds with $\gamma =1$ and $Q=F - G$.)
This optimality condition is readily interpreted: It states that $F$ is a 
fixed point of a projected gradient step with step size $\gamma$.

\paragraph{Optimality residual.} For any $Q \in \reals^{n \times m}$ with
$F=\Pi(Q)$, we define the (optimality) residual as
\[
r(Q) = \left\{ \begin{array}{ll} 
\min_{\gamma\geq 0}\|G-\gamma(F-Q)\|_F^2 & -FA^T \in \mathcal T \\
\infty & \mbox{otherwise},
\end{array} \right.
\]
where $\| \cdot \|_F^2$ denotes the squared Frobenius norm of a 
matrix, \ie, the sum of squares of its entries.
When $-FA^T \in\mathcal T $,
the righthand side is a quadratic function of $\gamma$, 
so the minimum is easily expressed explicitly as
\BEQ\label{e-res}
r(Q) =  \left\{ \begin{array}{ll} \|G\|_F^2
    -\frac{\Tr^2 G^T(F-Q)}{\|F-Q\|_F^2} 
&-FA^T \in \mathcal T,~ F \neq Q, ~\Tr G ^T(F-Q)\geq 0\\
\|G \|_F^2 & -FA^T \in \mathcal T,~ F = Q\text{ or }\Tr G ^T(F-Q)< 0\\
\infty & \mbox{otherwise}.
\end{array} \right.
\EEQ
Evidently $F=\Pi(Q)$ is optimal if and only if $r(Q)=0$.

\subsection{Related work}
\paragraph{Multicommodity network flow.}  
Historically, different forms of MCF problems have been formulated and studied. 
Starting from \cite{ford1958suggested} and \cite{hu1963multi} which studied a 
version with linear utility functions, which can be formulated as a linear 
program, later works develop nonlinear convex program 
formulations \cite{gautier1995forest,ouorou2000mean} and (nonconvex)
mixed integer program formulations
\cite{manfren2012multi,kabadurmus2016multi,zantuti2005wide} of MCF problems 
for different application purposes. These various forms of MCF have been widely used 
in transportation management \cite{erera2005global,mesquita2015rostering,
Rudi2016FreightTP}, 
energy and economic sectors \cite{Singh1978ADM,gautier1995forest,
manfren2012multi}, and network communication \cite{wagner2006fiber,kabadurmus2016multi,
layeb2017compact}. \cite{salimifard2022multicommodity} surveys over two hundred 
studies on MCF problems between 2000 and 2019. In this work, we focus on nonlinear 
convex formulation of MCF problems and develop GPU-compatible
algorithms for solving large problem instances. See \cite{ouorou2000survey} for a 
survey on nonlinear convex MCF problems. MCF models have very recently been exploited 
to design multi-GPU communication schedules for deep learning tasks 
\cite{Liu2024rethinking,basu2024efficientalltoallcollectivecommunication}, 
but the underlying MCF problems are solved with CPU-based solvers. 

\paragraph{First-order methods for convex optimization.}  First-order methods such 
as gradient descent algorithm, proximal point algorithm, primal-dual hybrid 
gradient algorithm, and their accelerated versions have been
exploited to tackle different forms of convex optimization problems. Compared 
to second-order methods which exploit Hessian information, first-order methods 
are known for their low computational complexity and are thus attractive for 
solving large-scale optimization problems. Recently, primal-dual hybrid gradient 
algorithm has been explored for 
solving large linear programs \cite{applegate2022practical,
lu2024cupdlpjlgpuimplementationrestarted,lu2024mpax} and optimal 
transport problems \cite{ryu2018vectormatrixoptimalmass} on GPUs. Other first-order 
methods such as ADMM have been exploited for designing 
GPU-accelerated optimizers for optimal power flow problems 
\cite{degleris2024gpuacceleratedsecurityconstrained,ryu2025gpu}. 

\paragraph{GPU-accelerated network flow optimization.} Specialized to GPU-based 
optimizers for  network flow optimization, \cite{wang2018toward} 
considers implementing a parallel routing algorithm on GPUs for SDN networks, 
which solves the Lagrangian relaxation of a mixed integer linear program. 
\cite{kikuta2015effective} implements a genetic method on GPUs for 
solving an integer linear program formulation of routing problem. \cite{zhang2023self} 
considers a linear program formulation of multicommodity network flow problems and 
constructs a deep learning model for generating new columns in delayed 
column generation method. \cite{WU201255} implements an asynchronous 
push-relabel algorithm for single commodity maximum network flow problem, 
which is CPU-GPU hybrid. \cite{liu2024keep} exploits exactly the same flow 
aggregation formulation of MCF following \cite{yin2019networkoptimizationunifiedpacket} 
as we do and trains a neural network model for minimizing 
unconstrained Lagrangian relaxation objective, and feeds the result 
as warm start to Gurobi \cite{gurobi} to get the final answer. 
\cite{Yarmoshik2024} integrates a source-based flow aggregation 
formulation of the multicommodity flow problem into 
solving the combined transportation 
model and exploits an accelerated variant of proximal alternating predictor-corrector 
algorithm. The authors claim that the proposed algorithm is GPU-friendly, but the 
numerical experiments are CPU-based, and involve small size networks. 
\cite{kubentayeva2023primaldualgradientmethodssearching} adopts a primal-dual 
gradient method for solving combined traffic models, which however is not 
GPU-oriented.

\subsection{Contribution}
Motivated by recent advancement of GPU optimizers, in this work we seek to 
accelerate large-scale nonlinear convex MCF problem solving with GPUs. 
Specifically, we adopt 
the MCF problem formulation in \cite{yin2019networkoptimizationunifiedpacket} 
(also described above) which is compactly matrix-represented and requires
fewer optimization variables by exploiting flow aggregation. We show that 
this specific problem formulation can be efficiently solved with first-order 
primal-dual hybrid gradient method when run on GPUs.

To the best of our knowledge, our work is the first to tackle exactly solving 
convex MCF problems on GPUs. Classic works for solving such 
large-scale MCF problems usually adopt Lagrangian relaxation for 
the coupling constraint and solve the resulting subproblems with smaller 
sizes in parallel (see, \eg, \cite{ouorou2000survey}). In our work, we 
do not exploit any explicit problem decomposition strategy and our 
algorithmic acceleration is mainly empirical and  depends on highly-optimized 
CUDA kernels for matrix operations. Moreover, we achieve problem size reduction via 
flow aggregation. Therefore our method has a simpler form which does not 
involve massive subproblem solving and synchronizing, and is also exact.

\subsection{Outline}
We describe our algorithm in \S \ref{s-pdhg}.
Experimental results, using our PyTorch implementation,
are presented and discussed in \S \ref{s-experiment}; very similar 
results obtained with our JAX implementation are given in appendix~\ref{s-jax}.
We conclude our work in \S \ref{s-conclusion}.  
The code, and all data needed to reproduce the results reported in this paper,
can be accessed at 
\begin{quote}
\url{https://github.com/cvxgrp/pdmcf}.
\end{quote}

\section{Primal-dual hybrid gradient} \label{s-pdhg}

\subsection{Primal-dual saddle point formulation}
We first derive a primal-dual saddle point formulation of the MCF 
problem \eqref{e-prob}. Let $\mathcal I$ denote the indicator 
function of $\mathcal F$, \ie, $\mathcal I(F)=0$ for $F \in \mathcal F$ and
$\mathcal I(F)=\infty$ otherwise.
We switch to minimizing $-U$ in \eqref{e-prob} to obtain the 
equivalent problem
\BEQ\label{e-prob-FT}
\begin{array}{ll}
\mbox{minimize} & -U(T)+\mathcal I(F) \\
\mbox{subject to} &  T=-FA^T, 
\end{array}
\EEQ
with variables $T$ and $F$.
We introduce a dual variable $Y \in \reals^{n \times n}$ 
associated with the (matrix) equality constraint. 
The Lagrangian is then
\[
\mathcal L(T,F;Y) = -U(T)+\mathcal I (F) -\Tr Y^T(T+FA^T)
\]
(see \cite[Chap.~5]{boyd2004convex}).
The Lagrangian $\mathcal L$ is convex in the primal variables 
$(T,F)$ and affine (and therefore concave) in the dual variable $Y$. 
If $(T,F;Y)$ is a saddle point of $\mathcal L$,
then $(T,F)$ is a solution to  problem \eqref{e-prob-FT} 
(and $F$ is a solution to the MCF problem \eqref{e-compact-prob}); the 
converse also holds. 

We can analytically minimize $\mathcal L$ 
over $T$ to obtain the reduced Lagrangian
\BEQ\label{e-L}
\hat {\mathcal L}(F;Y) = \inf_T \mathcal L(T,F;Y)=
-(-U)^*(Y) + \mathcal I(F) - \Tr Y^TFA^T,
\EEQ
where $U^*$ is the conjugate function of $U$ \cite[\S 3.3]{boyd2004convex}.
This reduced Lagrangian is convex in the primal variable $F$ and 
concave in the dual variable $Y$.
If $(F;Y)$ is a saddle point of $\hat {\mathcal L}$, then $F$ is a 
solution to the MCF problem \eqref{e-compact-prob} (see 
\cite[\S~1]{malitsky2018first}).
We observe that $\hat {\mathcal L}$ is convex-concave, with a 
bilinear coupling term.

\subsection{Basic PDHG method}
The primal-dual hybrid gradient (PDHG) algorithm, as first introduced 
in \cite{zhu2008efficient} and later studied in \cite{chambolle2011first,
chambolle2016ergodic}, is a first-order method for finding
a saddle point of a convex-concave function with bilinear 
coupling term.  The algorithm was extended to include over-relaxation 
in \cite[\S 4.1]{chambolle2016ergodic}, which has been observed to 
improve convergence in practice.
For \eqref{e-L}, PDHG has the form
\BEQ\label{e-pdhg-step}
\begin{array}{lll}
\hat F^{k+1/2} &=& \prox_{\alpha \mathcal I} (F^{k-1/2}+ \alpha Y^kA)\\
F^{k+1} &=& 2\hat F^{k+1/2}-F^{k-1/2}\\
\hat Y^{k+1} &=& \prox_{\beta (-U)^*} (Y^k - \beta F^{k+1}A^T)\\
F^{k+1/2} &=& \rho \hat F^{k+1/2} + (1-\rho) F^{k-1/2} \\
Y^{k+1} &=& \rho \hat Y^{k+1} + (1-\rho) Y^k
\end{array}
\EEQ
where $\prox_f(v)=\argmin_x(f(x)+(1/2)
\|x-v\|_2^2)$ denotes the proximal operator of $f$
\cite{parikh2014proximal}, 
$\alpha,\beta>0$ are positive step sizes satisfying 
$\alpha \beta \leq 1/\|A\|_2^2$,
and $\rho \in (0,2)$ is the over-relaxation parameter.

Reasonable choices for the parameters are
\[
\alpha = \beta = 1/\|A\|_2, \qquad \rho =1.9.
\]
(An upper bound on $\|A\|_2$ can be used in place of $\|A\|_2$.)


\paragraph{Convergence.}
In \cite{chambolle2016ergodic} it has been shown that when
there exists a saddle point of $\hat {\mathcal L}$,
$(F^k;Y^k)$ converges to a saddle point of $\hat {\mathcal L}$ 
as $k\rightarrow\infty$. For MCF the existence of an optimal
flow matrix and dual variable is known, 
so $F^k$ converges to an optimal flow matrix.
It follows that $r(F^{k-1/2}+\alpha Y^kA)\to 0$ as $k \to \infty$.
We note that $-FA^T \in \mathcal T$ only holds eventually.

\subsection{Proximal operators}
Here we take a closer look at the two proximal operators appearing 
in PDHG.

\paragraph{First proximal operator.}
We note that $\prox_{\alpha \mathcal I}$ appearing in the 
$\hat F^{k+1/2}$ update of 
\eqref{e-pdhg-step} is projection onto $\mathcal F$,
\[
\prox_{\alpha \mathcal I} (F) =\Pi (F).
\]
Since the constraints that define $\mathcal F$ separate 
across the columns of $F$, we can compute $\Pi(F)$ by
projecting each column $f_\ell$ of $F$
onto the scaled simplex $\mathcal S_\ell =
\{f \mid f \geq 0,~\ones^Tf \leq c_\ell\}$. 
This projection has the form
\[
\Pi_{\mathcal S_\ell}(f_\ell) = (f_\ell-\mu_\ell\ones)_+,
\]
where $\mu_\ell$ is the optimal Lagrange multiplier and
$( a )_+ = \max\{a,0\}$, which is applied elementwise to a vector.
The optimal $\mu_\ell$ is the smallest nonnegative value for which
$(f_\ell-\mu_\ell\ones)_+^T\ones \leq c_\ell$. 
This is readily found by a bisection algorithm; see
\S \ref{s-impl}.

\paragraph{Second proximal operator.}
The proximal operator appearing in the $\hat Y^{k+1}$ update step 
in \eqref{e-pdhg-step} can be decomposed entrywise, since $\beta (-U)^*$
is a sum of functions of different variables.
(The diagonal entries $-u_{ii}$ are zero, so
$(-\beta u_{ii})^*$ is the indicator function of $\{0\}$, and its 
proximal operator is the zero function.)
For each off-diagonal entry $i \neq j$ we need to evaluate
\[
\prox_{\beta (-u_{ij})^*} (y).
\]

These one-dimensional proximal operators are readily computed
in the general case.
For the weighted log utility $u(s)=w\log s$, we have
\[
\prox_{\beta(-u)^\ast}(y)=\frac{y-\sqrt{y^2+4\beta w}}{2}.
\]
For the weighted power utility $u(s)=ws^\gamma$,
$\prox_{\beta(-u)^\ast}(y)$ is the unique negative number $z$ for which
\[
(-z)^{c_1+2} +y(-z)^{c_1+1}-c_1c_2=0,
\]
where
\[
c_1=\frac{\gamma}{1-\gamma} >0,\quad c_2=\beta\left(
    \frac{1}{\gamma}-1\right)(w\gamma)^{\frac{1}{1-\gamma}}>0.
\]




\subsection{Adaptive step sizes}
In the basic PDHG algorithm \eqref{e-pdhg-step}, the step sizes 
$\alpha$ and $\beta$ are fixed.  It has been observed that 
varying them adaptively as the algorithm runs can improve 
practical convergence substantially \cite{applegate2022practical}.  
We describe our implementation of adaptive step sizes here.

We express the step sizes as
\[
\alpha^k =\eta/\omega^k, \qquad  \beta^k=\eta \omega^k,
\]
where $\eta \leq 1/\|A\|_2$ and
$\omega^k >0$ gives the primal weight.
With $\omega^k=1$ we obtain basic PDHG \eqref{e-pdhg-step}.

The primal weight $\omega^k$ is initialized as $\omega^0=1$ 
and adapted following \cite[\S 3.3]{applegate2022practical} as
\BEQ\label{e-primal_weight}
\omega^{k+1}=\left(\frac{\Delta_Y^{k+1}}
            {\Delta_F^{k+1}}\right)^\theta
\left(\omega^k\right)^{1-\theta},
\EEQ
where $\Delta_F^{k+1}=\|F^{k+1/2}-F^{k-1/2}\|_F$, $\Delta_Y^{k+1}=
\|Y^{k+1}-Y^k\|_F$ and $\theta$ is a parameter fixed as $0.5$ in 
our implementation.
The intuition behind the primal weight update 
\eqref{e-primal_weight} is to balance the primal and dual 
residuals;
see \cite[\S 3.3]{applegate2022practical} for details. 
In
\cite{applegate2022practical} the authors update $\omega$ 
each restart. We do not use restarts, and have found that updating
$\omega^k$ every $k^\text{adapt}$ iterations, when both 
$\Delta_F^{k}>10^{-5}$ and $\Delta_Y^{k}>10^{-5}$ hold,
works well in practice for MCF.  In our experiments we use $k^\text{adapt}=100$.
We can also stop adapting $\omega^k$ after some number of iterations,
keeping it constant in future iterations.  At least technically this 
implies that the convergence proof for constant $\omega$ holds for
the adaptive algorithm.

\paragraph{A simple bound on $\|A\|_2$.}
We can readily compute a simple upper bound on
\[
\|A\|_2=\sqrt{\lambda_{\max}(AA^T)},
\]
where $\lambda_{\max}$ denotes the maximum eigenvalue.
We observe that $AA^T$ is the Laplacian matrix associated with
the network, for which the well-known bound
\[
\lambda_{\max}(AA^T)\leq 2d_{\max}
\]
holds, where $d_{\max}$ is the largest node degree in the graph.
(For completeness we derive this in appendix \ref{s-laplacian-approx}.)
Thus we can take
\BEQ\label{e-eta}
\eta = 1/\sqrt{2d_{\max}}.
\EEQ

\subsection{Algorithm}\label{alg_sec}
We summarize our final algorithm, which we call PDMCF.
We set $r^0=+\infty$, $\alpha^0=\eta/\omega^0$, 
 and $\beta^0=\eta\omega^0$,
where $\eta$ is given in \eqref{e-eta} and  $\omega^0=1$.

\begin{algdesc}{\sc PDMCF} \label{pdhg-alg}
    {\footnotesize
    \begin{tabbing}
    {\bf given} $F^{-1/2},Y^0$, 
    parameter $\epsilon>0$.\\*[\smallskipamount]
    {\bf for $k=0,1,\ldots$} \\
    \qquad \= 1.\ \emph{Check stopping criterion.} Quit and return  
    $\hat F^{k-1/2}$
    if~$r^k<nm\epsilon$ 
    holds.\\
    \> 2.\ \emph{Basic PDHG updates \eqref{e-pdhg-step}.} \\ 
\>\qquad $\hat F^{k+1/2} =\Pi (F^{k-1/2}+\alpha^k Y^kA)$.\\
\>\qquad $F^{k+1} = 2\hat F^{k+1/2}-F^{k-1/2}$.\\
\>\qquad $\hat Y_{ij}^{k+1} = \left\{ \begin{array}{ll} \prox_{\beta^k(
    -u_{ij})^\ast}(Y_{ij}^k-
\beta^k(F^{k+1}A^T)_{ij})& j\neq i\\
0 & j=i. 
\end{array} \right. $ \\
\>\qquad $F^{k+1/2} = \rho \hat F^{k+1/2} + (1-\rho) F^{k-1/2}$. \\
\>\qquad $Y^{k+1} = \rho \hat Y^{k+1} + (1-\rho) Y^k$. \\

    \> 3. \emph{Adaptive step size updates \eqref{e-primal_weight}} 
    (if $k$ is multiple of $k^{\text{adapt}}$ 
    and $\Delta_F^{k+1},\Delta_Y^{k+1}>\tau$).\\
    \>\qquad $\omega^{k+1}=\left(\Delta_Y^{k+1}
            /\Delta_F^{k+1}\right)^\theta
\left(\omega^k\right)^{1-\theta}.$ \\
    \>\qquad $\alpha^{k+1}=\eta/\omega^{k+1},$ \qquad $\beta^{k+1}=
    \eta\omega^{k+1}.$
    \end{tabbing}}  
\end{algdesc}

\paragraph{Initialization.} We always take $F^{-1/2}=0$ and 
$Y^0=I-\ones\ones^T$. We can alternatively use 
a better guess of $F^{-1/2}$ and $Y$,
for example in a warm start, when
we have already solved a problem with similar data. We illustrate more on 
this in \S \ref{main_result}.

\paragraph{Stopping criterion.} Since $\hat F^{k+1/2}$ 
is result of projection onto $\mathcal F$, our optimality 
residual \eqref{e-res} has the form
\[
r^{k+1}=r(F^{k-1/2}+\alpha^k Y^kA).
\]
We consider the stopping criterion $r^k < nm\epsilon$,
\ie, the entrywise normalized residual 
$r^k/nm$ is smaller than a user-specified threshold $\epsilon$. 

\subsection{Implementation details}\label{s-impl}
\paragraph{Incidence matrix indexing.}
We only store the indices of the non-zero entries of $A$. 
Matrix multiplication with $A$ and 
$A^T$ can be efficiently computed by exploiting scatter
and gather functions, which are highly optimized 
CUDA kernels and are available in most major GPU computing languages.

\paragraph{Projection onto scaled simplex.}
To compute $\mu_\ell$ when $(f_\ell)_+^T\ones> c_\ell$, we follow 
\cite{held1974validation} and first sort 
$f_\ell$ from largest entry to smallest entry to form $f_\ell'$. We then 
find the largest index $t$ such that $f_{\ell t}'-((\sum_{i=1}^{t} 
f_{\ell i}'-c_\ell)/t)>0.$ Finally we take $\mu_\ell=(\sum_{i=1}^t 
f_{\ell i}'-c_\ell)/t$. 

Some recent work develops a more efficient method 
to compute the projection onto the simplex set (\cite{duchi2008} and 
\cite{condat2016} for example);
we adopt the simpler algorithm described above for implementation simplicity.

\section{Experiments} \label{s-experiment}

We run all our experiments on a single H100 
GPU with 80 Gb of memory supported by 26 virtual CPU cores and 241 Gb of RAM.
The results given below are for our PyTorch implementation; similar results,
reported in appendix~\ref{s-jax}, are obtained with our JAX implementation.

\subsection{Examples}\label{main_result}
\paragraph{Data and parameters.}
We consider weighted log utilities of form  
$u_{ij}(T_{ij})=w_{ij}\log T_{ij}$. We take $\log w_{ij}$ to be uniform 
on [$\log 0.3,\log 3$]. 
For network topology, we first create $n$ two-dimensional data
points $\xi_i\in\reals^2$, each denoted by $(\xi_{ix},\xi_{iy})$ for 
$i=1,\ldots,n$. We take $\xi_{ix}$ and $\xi_{iy}$ uniform on $[0,1].$ 
Then we add both edges $(\xi_i,\xi_j)$ and $(\xi_j,\xi_i)$ when either 
$\xi_i$ is among the $q$-nearest neighbors of $\xi_j$ or vice versa. 
For each edge $\ell$, we impose edge capacity $c_\ell$
where we take $\log c_\ell$ to be uniform on [$\log 0.5,\log 5$]. 

We use stopping criterion threshold 
$\epsilon=0.01/(n(n-1))$ for small to medium size problems and $\epsilon=
0.03/(n(n-1))$ for large size problems. We compare to CPU-based commercial 
solver MOSEK, with default settings.  MOSEK 
is able to solve the problems to high accuracy; we have checked that for 
all problem instances, the normalized utility differences
between results of PDMCF and MOSEK are no more than around $0.01$.
The pairwise normalized (optimal) utilities range between around $1$ and $10$,
which means that PDMCF finds flows that are between 0.1\% and 1\% suboptimal compared to
the flows found by MOSEK.
\paragraph{Small to medium size problems.} Table \ref{torch_small} shows 
runtime for both MOSEK and PDMCF 
required to solve  problem instances of various sizes. 
The column titled $nm$ gives the number of scalar optimization variables
in the problem instance.
We see that our implementation of PDMCF on a GPU 
gives a speedup over MOSEK of
$10\times$ to $1000\times$, with more significant 
speedup for larger problem instances. 
We also report runtime for PDMCF when run on CPU,
which is still quicker than MOSEK but 
with a significantly lower speedup.
Similar performance is also observed for our JAX implementation,
reported in appendix \ref{s-jax}. 

\begin{table}
\centering
\begin{tabular}{rrrr|rrrr} 
\multicolumn{4}{c}{problem sizes} & \multicolumn{3}{c}{timing (s)} & 
\multirow{2}{*}{iterations}  \\  
$n$ & $q$ & $m$ & $nm$ & MOSEK & PDMCF (CPU) & PDMCF (GPU)  \\
\hline
$100$ & $10$ & $1178$ & $1\times 10^5$  & $5$ & $1$ & $0.5$ & $490$  \\
$200$ & $10$ & $2316$ & $5\times 10^5$  & $23$ & $2$ & $0.7$ & $690$ \\
$300$ & $10$ & $3472$ & $1\times 10^6$  & $95$ & $6$ & $0.8$ & $840$ \\
$500$ & $10$ & $5738$ & $3\times 10^6$  & $340$ & $18$ & $1.1$ & $950$ \\
$500$ & $20$ & $11176$ & $6\times 10^6$  & $1977$ & $34$ & $1.4$ & $790$ \\ 
$1000$ & $10$ & $11424$ & $1\times 10^7$  & $2889$ & $1382$  & $19.5$ & $7220$  \\ 
$1000$ & $20$ & $22286$ & $2\times 10^7$  & $16765$ & $349$ & $5.1$ & $1040$  
\end{tabular}
\caption{Runtime table for small and medium size problems.}\label{torch_small}
\end{table}
\paragraph{Large size problems.} Table \ref{torch_large} shows runtime 
for several large problem instances. MOSEK 
fails to solve all these problems due to memory limitations.
PDMCF handles all these problem instances, with the largest one
involving $10^9$ variables.

\begin{table}
\centering
\begin{tabular}{rrrr|rrrr} 
\multicolumn{4}{c}{problem sizes} & \multicolumn{3}{c}{timing (s)} & 
\multirow{2}{*}{iterations}   \\  
$n$ & $q$ & $m$ & $nm$ & MOSEK & PDMCF (CPU) & PDMCF (GPU) \\
\hline
$3000$ & $10$ & $34424$ & $1\times 10^8$  & OOM & $7056$ & $96$ &   $4140$  \\
$5000$ & $10$ & $57338$ & $3\times 10^8$  & OOM & $19152$ & $395$ &  $3970$   \\
$10000$ & $10$ & $114054$ & $1\times 10^9$  & OOM & $87490$ &  $1908$& $4380$   
\end{tabular}
\caption{Runtime table for large size problems.}\label{torch_large}
\end{table}

\begin{figure}
    \begin{center}   \includegraphics[width=0.5\textwidth]{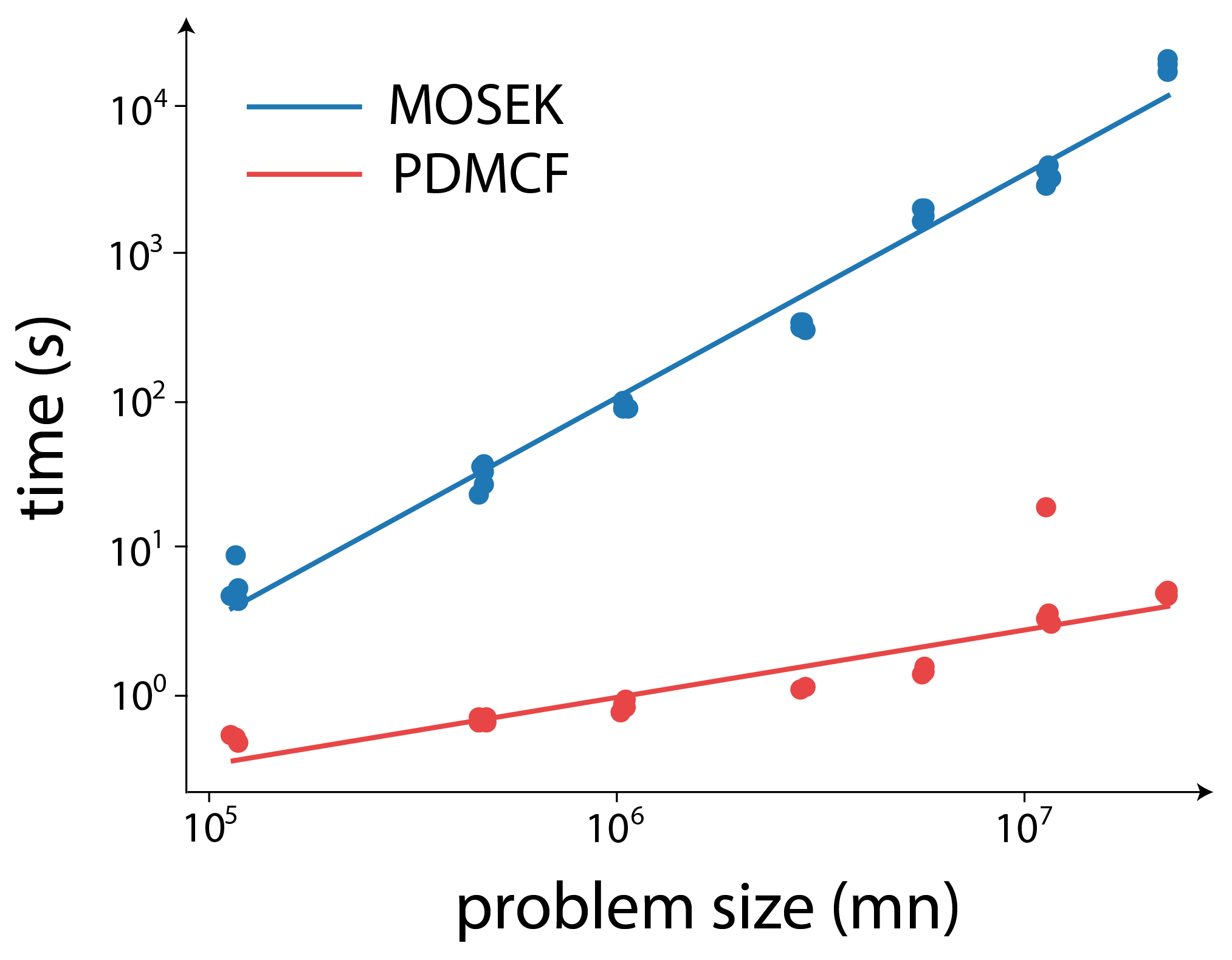}
    \end{center}
    \caption{Runtime plot for small and medium size problems.}\label{runtime_plot}
\end{figure}

\begin{figure}
    \begin{center}    
        \includegraphics[width=1.0\textwidth]{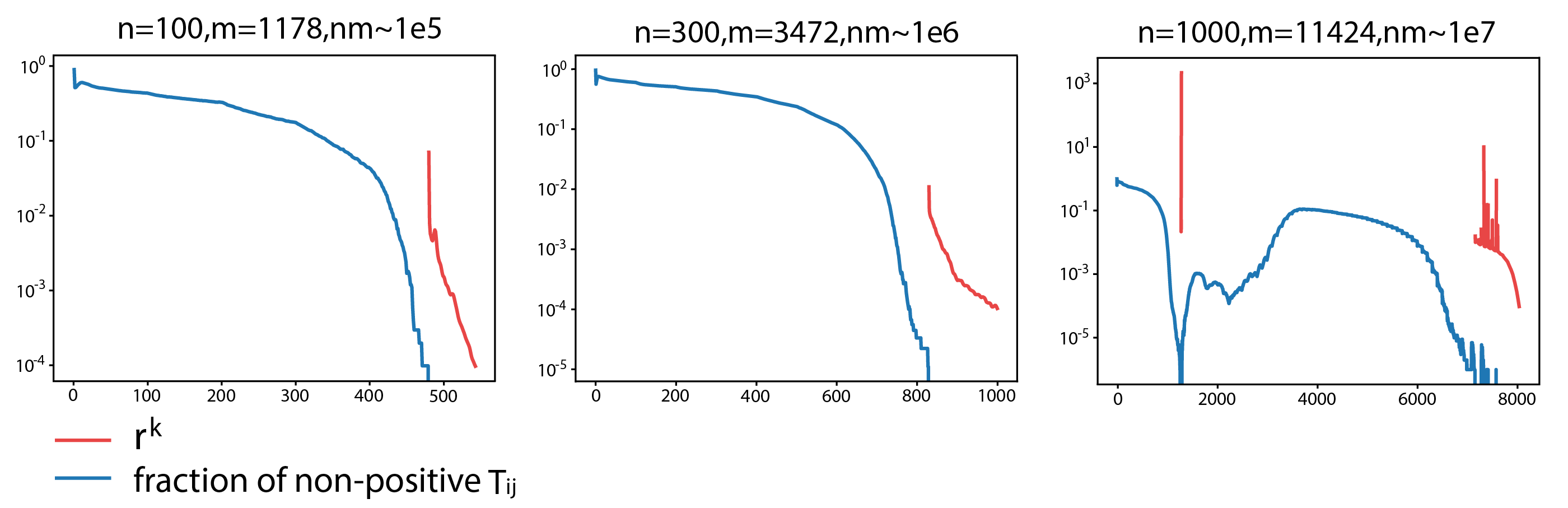}
    \end{center}
    \caption{Convergence plot for small and medium size problems.}
    \label{convergence_plot}
\end{figure}

\begin{figure}
    \begin{center}    
        \includegraphics[width=0.5\textwidth]{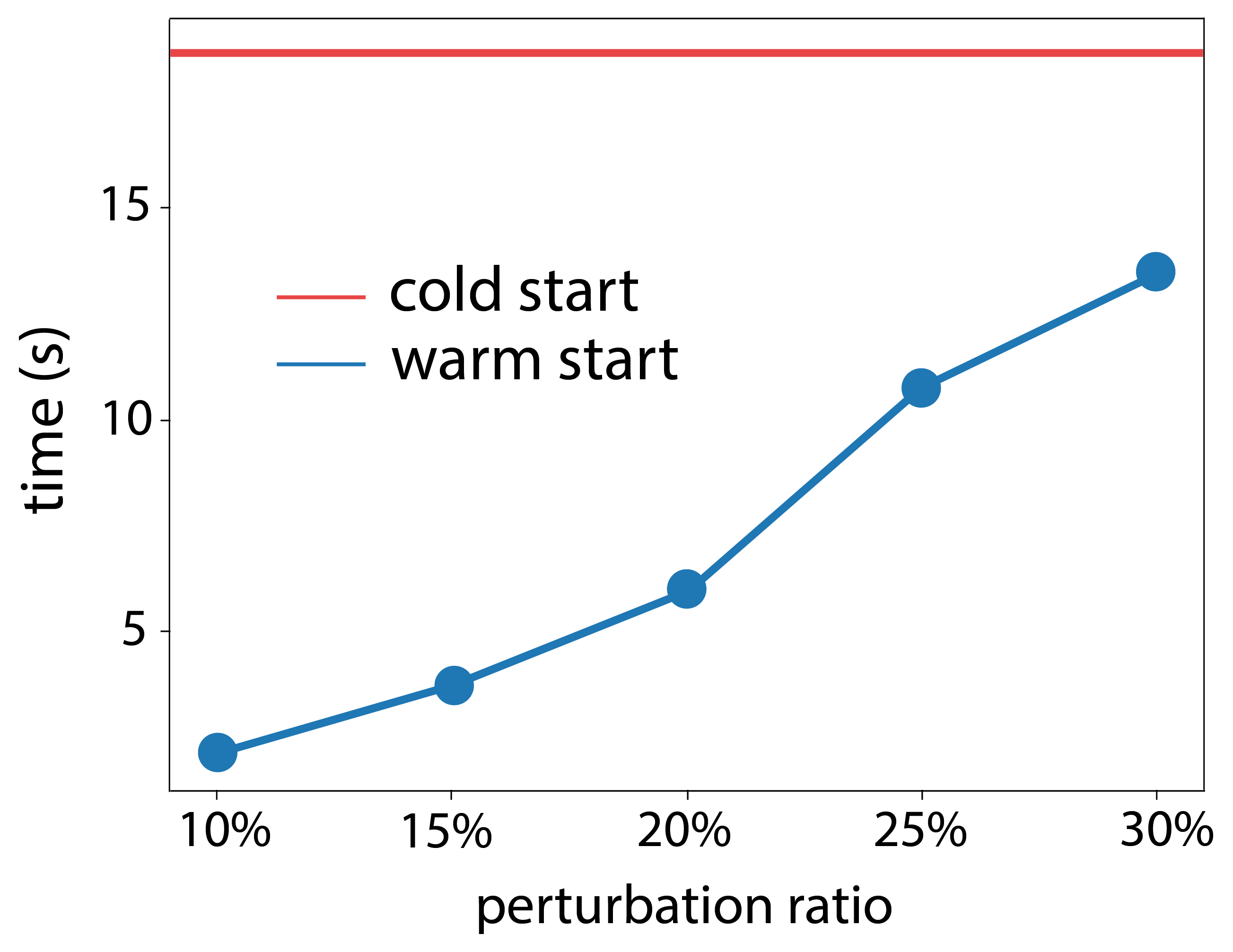}
    \end{center}
    \caption{Warm start plot for medium size problem.}\label{warm_plot}
\end{figure}

\paragraph{Scaling.} We scatter plot the runtime data for small
and medium problem instances in figure \ref{runtime_plot}.
Here we take $5$ problem instances generated by iterating over random 
seeds $\{0,1,2,3,4\}$ for the different $n,q$ values listed in 
table \ref{torch_small}. 
The $x$-axis represents optimization 
variable size $nm$ and the $y$-axis represents runtime in seconds. 
We plot on a log-log scale. 
The lines show the affine function fits to these data, with a 
slope around $1.5$ for MOSEK and around $0.5$ for PDMCF.
\paragraph{Convergence plot.} Figure \ref{convergence_plot} shows the 
convergence for three problem instances with variable 
sizes $10^5,10^6$, and $10^7$ with PDMCF, where the $x$-axis represents 
iteration numbers. 
Especially in the initial iterations we have infinite residual $r^k$
since $-FA^T \not\in \mathcal T$.  For those iterations we 
plot the fraction of nonpositive off-diagonal entries of $T$
in blue.
For feasible iterates we plot the (finite) residual, in red.
\paragraph{Warm start.} In \S\ref{alg_sec} we start with 
some simple initial $F^{-1/2}$ and $Y^0$. We  also test performance of 
PDMCF with warm starts.  In figure \ref{warm_plot} we present how runtime 
changes under different warm starts. To form these warm starts, for 
some perturbation ratio $\nu$, we randomly perturb 
entries of our utility weight matrix to derive $\tilde w_{ij}=(1\pm \nu)w_{ij}$,  
each with probability a half. We solve the multicommodity network flow 
problem with perturbed utility weight $\tilde w$ with PDMCF until we 
land at a feasible point 
$(F^{\text{feas}},Y^{\text{feas}})$ satisfying $(-F^{\text{feas}}A^T)_{ij}>0$ 
for all distinct $i,j$. We record the primal weight at this point as 
$\omega^{\text{feas}}$.  We then solve the desired multicommodity network flow 
problem with original 
utility weight $w$ with $F^{-1/2}=F^{\text{feas}},Y^0=Y^{\text{feas}}$ and 
$\omega^0=\omega^{\text{feas}}$. We note that setting $\omega^0=
\omega^{\text{feas}}$ is important for accelerated convergence, 
otherwise it usually requires similar number of iterations to converge as 
cold start if we simply set $\omega^0=1.$ In figure \ref{warm_plot}, we take 
problem instance with $n=1000,q=10$. $x$-axis stands for perturbation ratio 
$\nu$ and $y$-axis represents runtime in seconds. As can be observed, 
with perturbation ratio $\nu=10\%$, we harness $>80\%$ saving of runtime. 
Such savings keep decreasing to around $30\%$ when $\nu=30\%$, which makes 
sense given that larger perturbation indicates more different utility weights 
between original problem and perturbed problem, thus our warm start is expected 
to stay further from optimal solution to the original problem instance.
\section{Conclusion}\label{s-conclusion}
In this work, we present PDMCF algorithm which accelerates solving 
multicommodity network flow problems on GPUs. Our method starts with a 
destination-based formulation of multicommodity network flow problems 
which reduces optimization variable amount compared to classic problem 
formulation. We then apply PDHG algorithm to solve this destination-based 
problem formulation. Empirical results verify that our algorithm is GPU-friendly 
and brings up to three orders of magnitude of runtime acceleration compared 
to classic CPU-based commercial solvers. Moreover, our algorithm is 
able to solve ten times larger problems than those can be solved by commercial 
CPU-based solvers. 
\section*{Acknowledgements}
We thank Anthony Degleris and Parth Nobel for valuable discussions on 
implementation details. We also thank Demyan Yarmoshik for very useful 
feedback for the revision of our original manuscript.

\clearpage
\bibliography{references}
\clearpage
\appendix

\section{Upper bound on $\lambda_{\max}(AA^T)$}\label{s-laplacian-approx}
For a directed graph with incidence matrix $A$, $d_i=(AA^T)_{ii}$ 
is the degree of node $i$ and for $i \neq j$,
$-(AA^T)_{ij}$ is
the number of edges connecting node $i$ and node $j$, \ie, $2$ if both 
edges $(i,j)$ and $(j,i)$ exist. Note that  $\lambda_{\max}(AA^T)=
\max_{\|x\|_2=1}x^T(AA^T) x=\max_{x\neq 0}\frac{x^T(AA^T) x}{x^Tx}$.  
We have
\[
\begin{array}{ll}
x^T(AA^T)x &= \sum_i (AA^T)_{ii} x_i^2+\sum_{i\neq j} (AA^T)_{ij}x_ix_j\\
&=\sum_i d_i x_i^2+\sum_{i\neq j}(AA^T)_{ij}x_ix_j\\
&\leq \sum_i d_ix_i^2+\sum_{i\neq j} |(AA^T)_{ij}|(x_i^2/2+x_j^2/2)\\
&=\sum_i x_i^2(d_i+\sum_{j\neq i}|(AA^T)_{ij}|)\\
&=\sum_i 2d_i x_i^2\\
&\leq 2d_{\max}x^Tx.
\end{array}
\]
Therefore 
$\lambda_{\max}(AA^T)=\max_{x\neq 0}\frac{x^T(AA^T) x}{x^Tx}\leq 2d_{\max}$.
\section{JAX results}\label{s-jax}
The results shown in \S\ref{main_result} are for our
PyTorch implementation.  Here we provide the same results for 
our JAX implementation.
Tables \ref{jax_small} and \ref{jax_large} 
show the runtimes on the same  problem instances as reported in tables
\ref{torch_small} and \ref{torch_large}. We note 
that JAX's just-in-time (JIT) compilation adds runtime overhead 
for first-time function compilation and thus it does worse than 
its PyTorch counterpart on small size problems. 
The runtimes of these two versions 
are close for medium and large size problems, with JAX slightly 
slower.  

\begin{table}[ht!]
\centering
\begin{tabular}{rrrr|rrrr} 
\multicolumn{4}{c}{problem sizes} & \multicolumn{3}{c}{timing (s)} & 
\multirow{2}{*}{iterations}  \\  
$n$ & $q$ & $m$ & $nm$ & MOSEK & PDMCF (CPU) & PDMCF (GPU)  \\
\hline
$100$ & $10$ & $1178$ & $1\times 10^5$  & $5$ & $12$ & $5$ & $490$  \\
$200$ & $10$ & $2316$ & $5\times 10^5$  & $23$ & $57$ & $6$ & $690$ \\
$300$ & $10$ & $3472$ & $1\times 10^6$  & $95$ & $164$ & $6$ & $840$ \\
$500$ & $10$ & $5738$ & $3\times 10^6$  & $340$ & $548$ & $7$ & $950$ \\
$500$ & $20$ & $11176$ & $6\times 10^6$  & $1977$ & $890$ & $8$ & $790$ \\ 
$1000$ & $10$ & $11424$ & $1\times 10^7$  & $2889$ & $18554$ & $26$ & $7150$  \\ 
$1000$ & $20$ & $22286$ & $2\times 10^7$  & $16765$ & $5143$ & $15$ & $1040$  
\end{tabular}
\caption{Runtime table for small and medium size problems (JAX).}\label{jax_small}
\end{table}

\begin{table}
\centering
\begin{tabular}{rrrr|rrrr} 
\multicolumn{4}{c}{problem sizes} & \multicolumn{3}{c}{timing (s)} & 
\multirow{2}{*}{iterations}   \\  
$n$ & $q$ & $m$ & $nm$ & MOSEK & PDMCF (CPU) & PDMCF (GPU) \\
\hline
$3000$ & $10$ & $34424$ & $1\times 10^8$  & OOM & $106274$ & $139$ &   $4140$  \\
$5000$ & $10$ & $57338$ & $3\times 10^8$  & OOM & $382400$ & $421$ &  $3970$   \\
$10000$ & $10$ & $114054$ & $1\times 10^9$  & OOM & $1809517$ & $2078$ & $4380$   
\end{tabular}
\caption{Runtime table for large size problems (JAX).}\label{jax_large}
\end{table}
\end{document}